\input amstex
\input epsfx.tex
\frenchspacing
\documentstyle{amsppt}
\magnification=\magstep1
\baselineskip=14pt
\vsize=18.5cm
\footline{\hfill\sevenrm version 20071122} 
\def\cal{\Cal} 
\def\Tr{{\text{\rm Trace}}} 

\def\Pic{{\text{\rm Pic}}}

\def\End{{\text{\rm End}}}

\def\disc{{\text{\rm disc}}}

\def\mapright#1{\ \smash{\mathop{\longrightarrow}\limits^{#1}}\ }

\def\O{{\Cal O}}
\def\Que{\bold Q}
\def\CC{\bold C}
\def\FF{\bold F}
\def\HH{\bold H}
\def\Zee{\bold Z}
\def\Z{\bold Z}
\def\Q{\bold Q}
\def\H{\bold H}

\def\C{\bold C}
\def\P{\bold P}
\def\SL{{\text{\rm SL}}}
\def\Fq{\FF_q}
\def\Fp{\FF_p}
\def\FN{\FF_N}
\def\Fp{\FF_p}
\def\f{{\goth f}}
\def\gotha{{\goth a}}

\def\mapright#1{\ \smash{\mathop{\longrightarrow}\limits^{#1}}\ }

\newcount\refCount
\def\newref#1 {\advance\refCount by 1
\expandafter\edef\csname#1\endcsname{\the\refCount}}
\newref Abr  
\newref ALV  
\newref Be   
\newref BR   
\newref BSa  
\newref BSb  
\newref Cou  
\newref EM   
\newref ES   
\newref ESB  
\newref Gee  
\newref GS   
\newref HS   
\newref KSZ  
\newref KTK  
\newref LZ   
\newref LE   
\newref Mo   
\newref Po   
\newref Ruck 
\newref Sar  
\newref SSK  
\newref SC   
\newref ST   
\newref WE   
\topmatter
\title 
Constructing elliptic curves of prime order
\endtitle
\author Reinier Br\"oker, Peter Stevenhagen\endauthor
\address
University of Calgary, Department of Mathematics and Statistics,
2500 University Drive NW, Calgary, AB T2N 1N4, Canada
\endaddress
\address Mathematisch Instituut,
Universiteit Leiden, Postbus 9512, 2300 RA Leiden, The Netherlands
\email reinier\@math.ucalgary.ca, psh\@math.leidenuniv.nl\endemail
\endaddress
\thanks
This paper was submitte dfor publication on November 8, 2006.
It was completed at the Fields Institute in Toronto.
We thank this institute for its hospitality and support.
\endthanks
\abstract
We present a very efficient algorithm to construct
an elliptic curve $E$ and a finite field $\FF$ such that the order of
the point group $E({\FF})$ is a given prime number~$N$.
Heuristically, this algorithm only takes polynomial time
$\widetilde O((\log N)^3)$, and it is so fast that it may profitably
be used to tackle the related problem of finding elliptic curves with
point groups of prime order of prescribed size.

We also discuss the impact of the use of high level
modular functions to reduce the run time by large constant factors
and show that recent gonality bounds for modular curves imply limits
on the time reduction that can be obtained.
\endabstract
\endtopmatter

\document

\head 1. Introduction
\endhead

\noindent
For almost twenty years, the discrete logarithm problem in the group of
points on an elliptic curve over a finite field has been used as the basis
of elliptic curve cryptography.
Partly because of this application, the mathematically natural question of how
to generate elliptic curves over finite fields with a given number
of points has attracted considerable attention [\LZ, \KTK, \ALV, \BSa]. 
More in particular [\SSK, \KSZ], one is led to
the question of how to efficiently generate `cryptographic' elliptic curves
for which the order of the point group is a {\it prime\/} number.
For elliptic curves of prime order $N$, the discrete logarithm problem 
is currently supposed to be intractable for $N\gg 10^{60}$.

Section 2 deals with the problem of constructing a finite field $\FF$
and an elliptic curve $E/\FF$ having a prescribed prime number $N$ of
$\FF$-rational points.
We show that, on prime input $N$, such an elliptic curve 
can be constructed efficiently,
in heuristic polynomial time $\widetilde O((\log N)^3)$,
using traditional complex multiplication (CM) methods.
Here the $\widetilde O$-notation indicates that factors that are of
logarithmic order in the main term have been disregarded.
Note that $\widetilde O(X)$ for $X\to\infty$
is slightly more restrictive than
$O(X^{1+\varepsilon})$ for all $\varepsilon>0$.
The finite field $\FF$ over which $E$ is constructed will be of 
prime order $p$ for some $p$ sufficiently close to $N$.
The algorithm takes less time than algorithms that {\it prove\/}
the primality of the input $N$.
However, if the given input is known to be prime, the output
$E/\FF_p$ is guaranteed to be an elliptic curve over a prime field
$\FF_p$ having exactly $N$ points over $\FF_p$.
Because of its efficiency,
the range of our method amply exceeds the range of prime values
in current cryptographic use. 

In Section 3, we discuss the related problem of constructing an elliptic curve
that has a point group of prime order of prescribed {\it size\/}.
Unlike the earlier problem, this may be tackled efficiently
by `naive' methods that generate curves using trial and error
and exploit the efficiency of point counting on elliptic curves.
We describe the `traditional CM-algorithm' that constructs, on input of an integer
$k \in \Z_{\geq 3}$, an elliptic curve with prime order of $k$ decimal
digits, 
and show that the run time of this algorithm is $O(k^{4+\varepsilon})$
for every $\varepsilon>0$. 
It becomes $\widetilde O(k^4)$ if we are content with
probable primes instead of proven primes. As a consequence, we deduce that
the fastest way to tackle the problem in this Section
is to first {\it fix\/} a (probable)
prime $N$ of $k$ digits and then apply our CM-algorithm from Section 2
for that $N$.

{}From a practical point of view, CM-methods are hampered by the
enormous size of the auxiliary {\it class polynomials\/}
entering the construction, and since the time of Weber [\WE],
extensive use has been made of `small' modular functions
to perform CM-constructions.
We discuss the practical improvements of this nature in Section~4,
and show how recent results on the gonality of modular curves
imply upper bounds on the gain that can result from such methods.

A final section contains numerical illustrations of the methods
discussed.

\head 2. An efficient CM-construction
\endhead

\noindent
We start with the fundamental problem of realizing a prime number
$N>3$ as the group order of an elliptic curve $E$
defined over some finite field $\Fq$.
By Hasse's theorem, the order of the point group $E(\Fq)$ is an element of the
Hasse interval 
$$
{\cal H}_q = [q+1-2\sqrt{q},q+1+2\sqrt{q}] 
$$
around $q+1$.
The relation $N \in {\cal H}_q$ is actually symmetric in $N$ and $q$,
as we have $N \in {\cal H}_q \Longleftrightarrow q \in {\cal H}_N$. 
Consequently, a necessary condition for the
existence of a curve with $N$ points is that the Hasse interval ${\cal H}_N$
contains a prime power $q$.
As the set of integers $N$ for which ${\cal H}_N$ contains a non-prime
prime power $q$ is a zero density subset of $\Zee_{>0}$, we
may and will restrict to elliptic curves defined over {\it prime fields\/}
$\Fq = \Fp$.
If $p$ is a prime number in ${\cal H}_N$, then an elliptic
curve $E/\Fp$ with $\#E(\Fp) = N$ always exists. 
It follows from $p=N\in {\cal H}_N$ that
elliptic curves of prime order $N$ exist for every prime $N$,
but our algorithm will typically construct
curves over prime fields different from $\FN$.
This is certainly desirable from a cryptographic point of view,
as curves of order $N$ over $\FN$ are cryptographically unsafe:
the discrete logarithm problem on them can be transformed [\Ruck]
into a discrete logarithm problem for the {\it additive\/} group 
of $\FN$ that is easily solved.

Let $p$ be any prime in ${\cal H}_N$, and write $N=p+1-t$.
Then we have $t\ne0$, as the primes $p$ and $N>3$ are not
consecutive numbers.
It is well known that a curve $E/\Fp$ has
$N$ points over $\Fp$ if and only if the Frobenius morphism $F_p: E \rightarrow
E$ satisfies the quadratic equation
$$
F_p^2 - t F_p + p = 0 
$$
in the endomorphism ring $\End(E)$.
This means that the subring $\Z[F_p] \subseteq \End(E)$ generated by
Frobenius is isomorphic to the imaginary quadratic order $\O_\Delta$ of 
discriminant $\Delta = t^2-4p<0$, with $F_p$ corresponding to the
element $(t+\sqrt\Delta)/2 \in \O_\Delta$ of trace $t$ and norm $p$.
As $t$ is nonzero, the curve is ordinary.
Conversely, if the endomorphism ring $\End(E)$ of an ordinary
elliptic curve $E/\Fp$ contains an element $F$ of degree $p$
and trace $F+\hat F=t$, and therefore a subring isomorphic to
$\O_\Delta$, then one of the twists of $E$ over $\Fp$ has $N$ points.
Thus, constructing an elliptic curve having $N$ points over $\Fp$ is
the same problem as constructing an ordinary elliptic
curve over $\Fp$ for which the endomorphism ring is isomorphic to some
quadratic order containing $\O_\Delta$.

Over the complex numbers, the $j$-invariants of curves with endomorphism ring
isomorphic to $\O_\Delta$ are the roots of the {\it Hilbert class polynomial\/}
$$
P_\Delta = \prod_{[Q] \in \Pic(\O_\Delta)}(X - j(\tau_Q)) \in \Z[X].
$$
Here $j : \H \rightarrow \C$ is the classical modular function 
on the complex upper half plane $\H$ with Fourier
expansion $j(z) = 1/q+744 +\ldots$ in $q = \exp(2\pi i z)$, and the 
points $\tau_Q={-b+\sqrt\Delta\over 2a}\in \H$ correspond in the standard
way to the ideal classes
$$
[Q]=[\Zee\cdot a+\Zee\cdot {-b+\sqrt\Delta\over 2}]\in \Pic(\O_\Delta).
$$
The polynomial $P_\Delta$ has integer coefficients,
so it can be computed by approximating the 
roots $j(\tau_Q)\in\C$ with sufficient accuracy.  Alternatively, one can 
use $p$-adic algorithms [\Cou, \BR, \BSa] to compute $P_\Delta$.

The polynomial $P_\Delta$ splits completely modulo $p$, and its roots in
$\Fp$ are the $j$-invariants of the elliptic curves $E/\Fp$
with endomorphism ring isomorphic to $\O_\Delta$.
If $j_0\not = 0,1728 \in\Fp$ is one of these roots, then the elliptic
curve
$$
E : Y^2 = X^3 + aX - a\qquad\text{with }a = {27j_0\over 4(1728-j_0)}\in\Fp
\leqno{(2.1)}
$$
has $j$-invariant $j_0$.
If we have $N \cdot P = 0$ for our prime number $N$ and
$P = (1,1) \in E(\Fp)$, then $E(\Fp)$ has order $N$.
Otherwise the quadratic twist $E': Y^2 = X^3 + g^2aX - g^3a$ with
$g \in \Fp^*$ a non-square has $N$ points over $\Fp$.
In the special cases $j_0 = 0,1728$ there are a few more twists to consider.

As we only need $\End(E)$ to {\it contain\/} an order isomorphic to
$\O_\Delta$, we can replace $\Delta$ in the argument above by the
field discriminant $D = \disc(\Q(\sqrt{\Delta}))$.
For most~$t$, the discriminant $\Delta = \Delta(p) = t^2-4p$
is of roughly the same size as $p$ and $N$.
Moreover, the associated field discriminant $D$, which is essentially
the squarefree part of $\Delta$, will be of the same size as $\Delta$ itself
for most $\Delta$.
As computing the Hilbert class polynomial $P_D \in \Z[X]$,
which has degree $h(D)\approx \sqrt {|D|}$ and coefficients
of size $\widetilde O(\sqrt{|D|})$, takes time at least
linear in $D$, the CM-algorithm will have {\it exponential\/} run time
$\widetilde O(N)$ for `most' choices of primes $p\in {\cal H}_N$.

There is however a way to select primes $p\in {\cal H}_N$ for which the
field discriminant $D = D(p) = \disc(\Q(\sqrt{\Delta(p)}))$ is only of 
{\it polynomial size\/} in $\log N$.
What we want is
a discriminant $D$ such that the order $\O_D$ contains an element $\pi$ of
prime norm $p$ for which we have $N = p+1 - \Tr(\pi)=\text{\rm Norm}(1-\pi)$.
Exploiting the symmetry in $p$ and $N$ and writing $\alpha=1-\pi$,
we can also say equivalently that we want
an order $\O_D$ containing an element $\alpha$ of norm $N$ with the
property that $p = N+1-\Tr(\alpha)=\text{\rm Norm}(1-\alpha)$ is prime.
Note that if $\pi\in\O_D$ has prime norm $p>2$, then $\alpha=1-\pi$ will
have {\it even\/} norm in case the residue class field of the primes over 2
in $\O_D$ is the field of 2 elements.
For prime values $N>5$, or more generally for odd $N>5$,
this means that we can
only use discriminants $D$ congruent to $5$ modulo $8$.

In principle, one can find the {\it smallest\/} $D$ for which $\O_D$ contains
an element $\alpha$ of norm $N$ such that $\text{\rm Norm}(1-\alpha)$ 
is prime.
To do so, one splits the prime $N$ in the imaginary quadratic orders $\O_D$
with ${D\overwithdelims() N}=1$ as $(N)=\gotha\bar\gotha$ for descending
values of $D= -3, -11, -19, \ldots$ congruent to $5\bmod 8$
until we find a value of $D$
such that $\gotha=\alpha\O_D$ is principal with generator $\alpha$
and $N+1\pm \Tr(\alpha)=\text{Norm}(1\pm\alpha)$ is prime.
Now assume the standard {\it heuristical\/} arguments that
the prime $\gotha\subset \O_D$ over $N$ will be principal with `probability'
$1/h(D)$ and that $\text{Norm}(1\pm\alpha)\approx N$ will be prime
with `probability' $1/\log N$.
Then it is shown in [\BSb, Theorem 4.1] that the expected value of the 
smallest suitable discriminant $D$ found in this way will be
$$
D = \widetilde O((\log N)^2).
$$
Moreover, as the principality of the ideal $\gotha\subset \O_D$ lying over $N$
can be tested effiently using the 1908 algorithm of Cornacchia [\SC],
we can expect to find this $D$ in time $O((\log N)^{4+\varepsilon})$.

Cornacchia's algorithm explicitly computes the positive 
integers $x,y$ that satisfy
$$
x^2 - Dy^2 = 4N
$$
in case such integers exist. 
For $D <-4$, such $x, y$ are uniquely determined by $N$.
If found, the element $\alpha=(x+\sqrt D)/2\in \O_D$ has norm $N$, and we
hope that one of the
elements $\text{Norm}(1\pm\alpha)=N+1\pm x$ is prime.
Cornacchia's algorithm consists of the computation of a square root $x_0\bmod N$
of $D \bmod N$ followed by what is basically the Euclidean algorithm
for $x_0$ and $N$. 
It takes probabilistic time $\widetilde O((\log N)^2)$ for each $D$.
Performing Cornacchia's algorithm for $D=-3,-11,\ldots$ up to a bound
of size $(\log N)^2$ takes time $\widetilde O((\log N)^4)$,
and this dominates the run time of the algorithm.
We will lower the heuristic run time to $\widetilde O((\log N)^3)$
by applying an idea attributed to J. Shallit in [\Mo] to speed up the algorithm.

We start from the observation that $N$ splits into principal primes in
$\O_D$ if and only if $N$ splits completely in the Hilbert class field $H_D$ of 
$\Q(\sqrt{D})$. If this is the case,
then $N$ also splits completely in the genus field $G_D \subseteq H_D$,
which is
obtained by adjoining the square roots of $p^* = (-1)^{(p-1)/2}p$ to
$\Q(\sqrt{D})$ for all odd prime divisors $p \mid D$.
We have ${p^* \overwithdelims() N} = 
{N \overwithdelims() p} = 1$ for all odd primes dividing $D$,
and we can save time if we do not try increasing values of $D$ until
we hit the smallest suitable~$D$, but rather construct a suitable
discriminant $D$ from a
generating set of `good' primes $p$ for which we know that $p^*$ is
a square modulo $N$.
If we only consider primes $p$ of size $O(\log N)$, the time
needed to compute the values $\sqrt{p^*} \bmod N$ for these primes
is $\widetilde O((\log N)^3)$. 

Our algorithm consists of multiple `search rounds' for a suitable
discriminant~$D$,
where in each round we increase the size of the `basis' of primes we use.
First we take the primes between $0$ and $\log N$ and see whether we can
find a suitable $D \equiv 5 \bmod 8$ with $|D| < (\log N)^2$ a product
of primes
from this basis. 
If no such $D$ exists, we add the `good' primes 
between $\log N$ and $2 \log N$ to our basis, and look for a
suitable $D$ with $|D| < (2 \log N)^2$
created from this enlarged basis, and so on.
In this way, we encounter in the $r$-th round all discriminants
$D$ with $|D| < (r \log N)^2$ that are products of
prime factors below $r\log N$.
Asymptotically (cf. the `analytic tidbit' in [\Po]), this is a {\it positive\/}
fraction $1-\log 2 \approx 0.30685$
of all discriminants below $(r \log N)^2$.
As the smoothness properties of $D$ play no role in
our heuristics,
we still expect to find a suitable discriminant of size $\widetilde
O((\log N)^2)$.
Thus, we expect the algorithm below to terminate after a number of
rounds that is polynomial in $\log\log N$.
In practice (cf. Section 5), this number is usually 1.
\medskip\noindent
{\bf 2.2. Algorithm.} 
\hfil\break
Input: a prime number $N$.\hfil\break
Output: a prime number $p$ and an elliptic curve $E/\Fp$ with
$\#E(\Fp) = N$.
\vskip5pt\parindent=.6cm
\item{\bf 1.}
Put $r \leftarrow 0$, and create an empty table $S$.
\item{\bf 2.}
Compute for all odd primes $p\in [r\log N, (r+1)\log N]$ that
satisfy ${N \overwithdelims() p} = 1$
a square root $\sqrt{p^*} \bmod N$,
and add the pairs $(p^*, \sqrt{p^*} \bmod N)$ to the table~$S$.
\item{\bf 3.}
For each product
$(D,\sqrt{D} \bmod N)=(\prod_i p_i^*, \prod_i \sqrt{p_i^*} \bmod N)$
of distinct elements of $S$ that satisfies $\prod_i p_i^* < (r\log N)^2$
and $D \equiv 5 \bmod 8$,
do the following.
\itemitem{\bf 3a.}
Use the value $\sqrt{D} \bmod N$ and Cornacchia's algorithm 
to compute $x, y>0$ satisfying $x^2-Dy^2 = 4N$.
\itemitem{\bf 3b.}
For each solution found in step 3a, test whether $p = N+1
\pm x$ is a probable prime.
If it is, compute the Hilbert class polynomial $P_D \in \Z[X]$,
compute a root $j_0$ of $\overline P_D \in \Fp[X]$, return the twist
of the elliptic curve (2.1) that has $N$ points, and stop.
If no root or no twist is found, then $p = N+1\pm x$ is not prime
and we continue with the next solution.
\item{\bf 4.} Put $r \leftarrow r+1$ and go back to step 2.
\medskip\noindent
A heuristic analysis of the Algorithm above leads to the following.
\proclaim
{2.3. Theorem}
On input of a prime $N$,
Algorithm 2.2 returns a prime $p$
and an elliptic curve $E/\Fp$ with $\#E(\Fp) = N$.
Under heuristic assumptions, its run time is $\widetilde O((\log N)^3)$.
\endproclaim
\noindent{\bf Proof. }
As the smoothness properties of $D$ are irrelevant in the
heuristic analysis detailed in [\BSb],
the smallest suitable $D$ found by our Algorithm,
which restricts to the positive density subset of discriminants,
will be of size $\widetilde O((\log N)^2)$.
The expected number $r$ of rounds of our Algorithm
will therefore be small, at most polynomial in $\log\log N$,
and in view of our $\widetilde O$-notation we may prove our Theorem
by focusing on the time needed for a single round of the Algorithm,
which consists of Steps 2 and 3.

In Step 2 we have to find primes up to $(r+1)\log N$. As we only need
to test primality of integers of size $\widetilde O(\log N)$,
the time needed to find these primes is negligible.
For all the $\widetilde O(\log N)$ primes we find, we need to test
which primes are `good', i.e., which primes satisfy
${N \overwithdelims() p} = 1$. The time needed for this computation is
also negligible. The bottleneck in Step 2 is the computation of the
square roots of $p^* \bmod N$ for the good primes $p$. Each square 
root computation takes time $\widetilde O((\log N)^2)$, so Step 2 takes 
time $\widetilde O((\log N)^3)$.

For each of the $O((\log N)^2)$ products $(D, \sqrt D\bmod N)$
formed in Step 3, we run the Euclidean algorithm part of
Cornacchia's algorithm in Step 3a in time 
$\widetilde O(\log N)$.
This takes time $\widetilde O((\log N)^3)$.
We expect to find $O(\log N)$ solutions $(x,y)$ from Step 3a
for which we have to test primality of $N+1\pm x$ in Step 3b.
A cheap Miller-Rabin test, which takes time $\widetilde O((\log N)^2)$,
suffices for our purposes, and leads to a total time 
$\widetilde O((\log N)^3)$ spent on primality testing.

Once we encounter a probable prime $p=N+1\pm x$ for some discriminant~$D$, 
we compute the Hilbert class polynomial $P_D$. As $D$ is of size
$O((r\log N)^2)$, this takes time $\widetilde O((\log N)^2)$. 
Computing a root $j_0$ of $P_D$, a polynomial of degree $h(D)=\widetilde O(\log N)$,
modulo the prime $p\approx N$ once more takes time $\widetilde O((\log N)^3)$.
To test which curve of $j$-invariant $j_0$ has $N$ points, we may have
to compute all isomorphism classes over $\Fp$ of elliptic curves
with $j$-invariant $j_0$ until we find one.
There are at most 6 of these classes (`twists'),
and for the class of $E$ we need to test the equality $N \cdot P=0$ for
a point $P$ on~$E$.
This only takes time $\widetilde O((\log N)^2)$, and we conclude that the entire
round of the algorithm runs in time $\widetilde O((\log N)^3)$.

Even though we have only found a {\it probable\/} prime $p$ in
the beginning of Step~3b, the equality $N \cdot P=0$ on $E$
tested in this Step exhibits a point of order $N$ on $E$, which
{\it proves\/} that $p$ is actually prime.
\hfill$\square$
\medskip\noindent
The low asymptotic running time of our Algorithm is
illustrated by the size of some of the examples in Section 5.
As several steps in the algorithm are no faster than
$\widetilde O((\log N)^3)$, it seems that we have obtained an optimal result
for a CM-solution to our problem.

\head 3. Point groups of given prime size
\endhead

\noindent
Closely related to the problem of constructing an elliptic curve of 
prescribed prime order $N$ is the problem of constructing a curve for
which the group order is a prime in a given interval.
For concreteness sake, we take the interval as $[10^{k-1},10^k)$, so the
problem becomes the efficient construction of
an elliptic curve over a finite field such that the group
order is a prime of exactly $k$ decimal digits.

If we insist on a curve with {\it proven\/} prime order, we cannot hope
for an algorithm with a faster run time than $O(k^4)$, since the fastest
known algorithm [\Be] to rigorously prove primality of an integer 
$N \approx 10^k$ has expected run time
$O((\log N)^{4+\varepsilon}) =  O(k^{4+\varepsilon})$ for all
$\varepsilon>0$.
The naive algorithm of selecting a prime $p$ of $k$ decimal digits and trying 
random elliptic curves over $\Fp$ until we find one of prime order 
already has a heuristic run time that comes close to this `optimal run time'.
Indeed, counting the number of points of an elliptic curve $E/\Fp$ takes
heuristic time $\widetilde O((\log p)^4)$ using the improvements made
by Atkin and Elkies to Schoof's original point counting algorithm [\SC]. 
Even though the distribution of group orders of elliptic curves over $\Fp$
over the Hasse interval is not exactly uniform, it follows as in
[\LE, Section 1] that, heuristically,
we have to try $O(\log p)$ curves over $\Fp$ until we find one of prime
order. This leads to a heuristic run time $\widetilde O(k^5)$.

As was noted by many people [\EM, \KSZ], we can also use complex multiplication
techniques to tackle the problem. 
Unlike our Algorithm 2.2, which starts with a desired prime value $N$
for the group order and computes a suitable prime field~$\Fp$ over which
the curve can be constructed, these algorithms compute primes~$p$
splitting into principal primes $\pi$ and $\bar\pi$ in some {\it fixed\/}
quadratic ring $\O_D$, and construct a curve over $\Fp$ having
CM by $\O_D$ and $N$ points when $N=\text{\rm Norm}(1-\pi)$ is found to be 
prime. As before, we can test whether a given prime $p$ splits into
principal primes in $\O_D$ by computing a value of $\sqrt{D} \bmod p$ for
${ D\overwithdelims() p } = 1$ and
applying Cornacchia's algorithm. In case $\O_D$ has class number $1$, i.e.,
for $D = -3, -11, -19, -43, -67, -163$, we can see whether $p$ splits
in $\O_D$ by only looking at $p$ mod $D$.

Subject to the congruence condition $D \equiv 5 \bmod 8$, we can take
{\it any\/} fundamental discriminant. The run time depends on the value
of $D$ we choose, the value $D=-3$ being `optimal'. For cryptographic
purposes we need to select $D$ such that the class number of $\O_D$ is
at least 200, cf.\ Section 5.\par
\ \par
\noindent
{\bf 3.1. Algorithm.} 
\hfil\break
\noindent
Input: an integer $k \in \Z_{\geq 3}$, and a negative discriminant $D \equiv
5 \bmod 8$.\par
\noindent
Output: primes $p,q$ of $k$ decimal digits and an elliptic curve $E/\Fp$ 
with CM by $\O_D$ and $\#E(\Fp) = q$.\par
\vskip5pt\parindent=.6cm
\item{\bf 1.} Compute $P_D \in \Z[X]$.
\item{\bf 2.}
Pick a random probable prime $p$ that splits into principal primes in $\O_D$
and satisfies
$$
10^{k-1}+ 2\cdot {10^{k-1 \over 2}}< p <  10^k-2\cdot {10^{k\over 2}}.
$$
\item{\bf 3.}
Write $p = \pi \overline \pi \in \O_D$. If $q=\text{\rm Norm}(1-
\varepsilon\pi)$ is a probable prime for some $\varepsilon\in \O_D^*$, prove
the primality of $q$, compute a root $j\in\Fp$ of $P_D \in \Fp[X]$ and
return an elliptic curve $E/\Fp$ with $j$-invariant $j$ with $q$ points. Else, go 
to Step 2.\par
\ \par\noindent
A heuristic analysis of the Algorithm above leads to the following.\par
\proclaim
{3.2. Theorem}
On input of an integer $k \in \Z_{\geq 3}$ and a negative
discriminant $D \equiv 5 \bmod 8$,
Algorithm 3.1 returns primes $p,q$ of $k$ decimal digits 
each and an elliptic curve $E/\Fp$ with CM by $\O_D$ and $\#E(\Fp) = q$. 
Under heuristic assumptions, the run time for fixed $D$ is
$O(k^{4+\varepsilon})$ for every $\varepsilon>0$.
\endproclaim
\noindent{\bf Proof. }
To prove that the output of Algorithm 3.1 is correct,
we only need to check that the norms $q$ found in Step 2 have $k$
decimal digits.
This follows from Hasse's theorem $q\in{\Cal H}_p$ and the choice
of our interval for $p$.

In Step 1 we have to find a prime $p$ of $k$ decimal digits that splits into
principal primes in $\O_D$.
Finding a probable prime $p$ of $k$ digits takes time 
$\widetilde O(k^3)$, and with positive probability $(2h(D))^{-1}$ such
a prime $p$ splits into principal primes in $\O_D$.
For each $p$ found we can test this in time $\widetilde O(k^2)$
by computing a value $\sqrt{D} \bmod p$ in case it exists,
and use it to apply the Euclidean algorithm part of Cornacchia's algorithm.
 
If $p$ factors as $p = \pi\overline\pi$ in $\O_D$, the `probability' that 
$\text{\rm Norm}(1-\varepsilon\pi)$ is prime is about $1/k$.
We expect that we need to perform Step 2 roughly $k$ times,
and except for the primality proof of $q$ this takes us time $\widetilde O(k^4)$.

A rigorous primality proof of $q$ in Step 3 takes time $O(k^{4 + \varepsilon})$
for every $\varepsilon>0$. Just as in Theorem 2.3, this also proves
the primality of $p$.  
\hfill$\square$
\medskip\noindent
The proof shows that if we only insist that $p,q$ are {\it probable\/} 
primes of $k$ digits, the run time becomes $\widetilde O(k^4)$. This is
{\it slower\/} than Algorithm 2.2. The fastest way of constructing a curve
for which the group order is a probable prime of $k$ digits is therefore to 
find a random probable prime $N$ of $k$ digits and then run Algorithm 2.2 
on this input. Indeed, finding a probable prime $N$ takes time 
$\widetilde O(k^3)$, and so does the application of Algorithm 2.2 on $N$.

\head 4. Class invariants and gonality
\endhead
\noindent
In large examples,
the practical performance of Algorithm 2.2 is hampered by the
computation of a Hilbert class polynomial $P_D$ in Step 3b.
As we noted already, the run time $\widetilde O(|D|)$ needed for computing $P_D$
cannot be seriously improved, as the degree $h(D)$ of $P_D$ is of order of
magnitude $\sqrt{|D|}$ by the Brauer-Siegel theorem, and the number of digits
of its coefficients has a similar order of magnitude $\sqrt{|D|}$.
However, already for the moderately small values of $D$ used by our algorithm,
the coefficients of $P_D$ are notoriously large.

It was discovered by Weber [\WE] that one can often work with `smaller' modular 
functions than the $j$-function to generate the Hilbert class field $H_D$
of $\Que(\sqrt D)$.
There are many of these functions, and each of them works for
some positive proportion of discriminants.
A good example is provided by the Weber function
$\f = \zeta_{48}^{-1} \eta({z+1\over 2})/ \eta(z)$,
which is related to $j$ by 
an irreducible polynomial relation 
$$
\Psi(\f, j)=(\f^{24} -16)^3 - j \f^{24}=0
$$
of degree 72 in $\f$ and degree 1 in $j$.
It can be used for all $D\equiv 1\bmod 8$ coprime to~3.
For $D=-71$, the value $\f(\tau)$ for an appropriate generator $\tau$
of $\O_{-71}=\Zee[\tau]$ has the irreducible polynomial
$$
P^\f_{-71}= X^7+X^6-X^5-X^4-X^3+X^2+2X+1 \in \Z[X]
$$
that requires less precision to compute from its complex zeroes than
it does to compute the Hilbert class polynomial
$$\eqalign{
P_{-71}&=  \;X^7 + 313645809715\;X^6 -
         3091990138604570\;X^5 \cr&\quad
        + 98394038810047812049302\;X^4
        - 823534263439730779968091389\;X^3 \cr &\quad
        + 5138800366453976780323726329446\;X^2 \cr&\quad
        - 425319473946139603274605151187659\; X \cr &\quad
        + 737707086760731113357714241006081263 \cr
        }
$$ 
coming from the $j$-function.
The polynomials $P_{-71}$ and $P^\f_{-71}$ have
the same type of splitting behavior modulo primes as they
generate the same field $H_{-71}$ over $\Que(\sqrt{-71})$.
Moreover, the zeroes modulo $p$ of $P^\f_{-71}$
readily give the zeroes of $P_{-71}$
modulo $p$ by the formula $j= \f^{-24}(\f^{24} -16)^3$.
A significant speed up in the practical performance of 
CM-algorithms can be obtained by using functions such as $\f$ instead of $j$.

In cases where the value $f(\tau)$ of a modular function $f$ at
some $\tau \in \Q(\sqrt{D})$ generates the Hilbert class field $H_D$
over $\Que(\sqrt D)$, we call $f(\tau)$ a {\it class invariant\/}. Class
invariants have been well studied, and it is now a rather mechanical process 
[\ST, \GS] to check for which $D$ class invariants can be obtained 
from a given modular function $f$, and, in case $f(\tau)$ is a 
class invariant for $\Que(\sqrt D)$, to find its Galois conjugates
and to compute its minimal polynomial $P^f_D$ over $\Q$. 

If $f$ yields class invariants, the logarithmic height of the 
zeroes of $P_D^f$ will asymptotically, for $D\to-\infty$, differ from
those of $P_D = P_D^j$ by some {\it constant factor\/} depending
on the function $f$.
This is the factor we gain in the size of the coefficients of $P_D^f$
when compared to $P_D$.
For the Weber function $\f$ above, we get class invariants for
discriminants $D \equiv 1 \bmod 8$ not divisible by 3, and the length
of the coefficients is a factor $72$ smaller for $P_D^\f$ than it is for $P_D$.
For other discriminants, such as the discriminants congruent to $5\bmod 8$ from
the previous sections, similar but somewhat smaller factors may be
gained by using double eta-quotients
$\eta(z/p)\eta(z/q)\eta(z)^{-1}\eta(z/pq)^{-1}$ as in [\ES].

%
The `reduction factor' that is obtained when using a modular function
$f$ instead of $j$ depends on the degree of the irreducible
polynomial relation $\Psi(j, f)=0$ that exists between $j$ and $f$.
In terms of the polynomial $\Psi(j, f)\in\CC[X]$, we define the
{\it reduction factor\/} of our modular function $f$ as
$$
r(f)=
{\deg_f(\Psi(f,j))\over \deg_j(\Psi(f,j))}.
$$
By [\HS, Proposition B.3.5], this is, asymptotically, the {\it inverse\/}
of the factor 
$$
\lim_{h(j(\tau)) \rightarrow \infty}
{h(f(\tau)) \over h(j(\tau))}.
$$
Here $h$ is the absolute logarithmic height, and we take the limit over
all CM-points $\SL_2(\Z)\cdot\tau\in\HH$, ordered by the absolute value of
the discriminant of the associated CM-order.
The reduction factor 72 obtained for the Weber function above is close
to optimal in view of the following theorem.
\proclaim
{4.1. Theorem}
The reduction factor of a modular function $f$ satisfies
$$
r(f)\le {800/7}\approx 114.28.
$$
If Selberg's eigenvalue conjecture in [\Sar] holds, then we have
$$
r(f)\le 96.
$$
\endproclaim
\noindent{\bf Proof. }
Let $f$ be modular of level $N \geq 1$, and $\Gamma(f)\subset \SL_2(\Z)$
the stabilizer of $f$ inside $\SL_2(\Z)$.
Then $\Gamma(f)$ contains the principal congruence
subgroup $\Gamma(N)$ of level~$N$, and the inclusions
$\Gamma(N)\cdot\{\pm1\}\subset \Gamma(f)\subset \SL_2(\Z)$
correspond to coverings 
$$X(N)\mapright{} X(f)\mapright{j} \P^1_\CC$$
of modular curves.
Here $X(N)$ is the full modular curve $X(N)$ of level $N$,
which maps to the $j$-line $\P^1_\CC$ under $j$.
This map factors via the intermediate modular curve $X(f)$,
which has function field $\CC(j, f)$.
The Galois theory for the function fields shows
that the degree of the map $j: X(f)\to \P^1_\CC$ is equal to
$$
[\SL_2(\Z):\Gamma(f)]=[\CC(j , f):\CC(j)]=\deg_f(\Psi(f, j)).
$$
We now consider the {\it gonality\/} $\gamma(X(f))$ of the modular curve
$X(f)$, i.e., the minimal degree of a non-constant morphism
$\pi: X(f)\to \P^1_\CC$.
Abramovich [\Abr] proved in 1996 that the gonality of {\it any\/}
modular curve $X_H$ corresponding to some congruence subgroup
$H\subset \SL_2(\Zee)$ is bounded from below by 
$c\cdot [\SL_2(\Z):H]$ for some universal constant $c>0$.
His proof yields the value $c=7/800$, and under assumption of 
Selberg's eigenvalue conjecture [\Sar] the constant $c$ can be taken
equal to $1/96$.

For our curve $X(f)$, the rational map $f: X(f)\to \P^1_\CC$ has degree
$$
[\CC(j , f):\CC(f)]=\deg_j(\Psi(f,j)),
$$
and this degree is at
least $\gamma(X(f))$.
We can now use Abramovich's lower bound to obtain
$$
r(f) = {\deg_f(\Psi(f,j))\over \deg_j(\Psi(f,j))}
\leq {[\SL_2(\Z):\Gamma(f)]\over \gamma(X(f))}\leq {1\over c}.
$$
The proven value $c=7/800$ and its conditional improvement
$c=1/96$ yield the two statements of our theorem.\hfill$\square$
\medskip\noindent
We do not know whether the value 96 is attained for some function $f$.
The factor 72 of Weber's function is the best we know.

\head 5. Numerical examples
\endhead
\noindent
We illustrate Algorithm 2.2 by constructing an elliptic curve having exactly
$$
N = 123456789012345678901234567890123456789012345678901234568197
$$
points. The integer $N \approx 10^{60}$ is prime, and the discrete logarithm
problem is believed to be hard for such a curve.

We have $\log N \approx 136$ and there are $15$ odd primes $p<136$ with
${ p^* \overwithdelims() N } = 1$. We compute and store $\sqrt{p^*} \bmod N$
for these primes, and we try to find a discriminant
$D\equiv5\bmod8$ built from primes out of this `basis' such that $N$ splits as $N=\alpha
\overline\alpha$ in the order $\O_D$ and such that $N+1\pm\Tr(\alpha)$ is prime.
For $D = -41\cdot 59= -2419$ we find a solution
\medskip
\centerline{\vbox{\halign{$#$&$\thinspace#$\cr
x = & 531376585512740287835890668303\cr
y = & 9349802208089011828618119329\cr}}}
\medskip\noindent
to the norm equation $x^2 - Dy^2 = 4N$ for which $p = N+1+x$ is prime. 

The class group $\Pic(\O_D)$ is cyclic of order $8$. The Hilbert class
polynomial $P_D$ has degree $8$, and coefficients of up to $119$ decimal
digits. It splits completely modulo 
$$
p = 123456789012345678901234567890654833374525085966737125236501,
$$
and any of its zeroes is the $j$-invariant of a curve having $N$ points.
With $a = 112507913528623610837613885503682230698868883572599681384335\in\Fp$,
the elliptic curve $E_a$ given by
$$
Y^2 = X^3+aX-a
$$
has $N$ points, as may be checked by computing $N \cdot (1,1) = 0 \in E_a(\Fp)$.

We can speed up the algorithm by computing a `smaller' polynomial
than the Hilbert class polynomial.
We are in the case where $3$ does not divide $D=-2419$, and here
the cube root $\gamma_2$ of the $j$-function can be shown to yield class invariants. 
The polynomial $P_{-2419}^{\gamma_2}\in\Z[X]$ has coefficients up to $40\approx 119/3$ decimal 
digits. For a root $x\in\Fp$ of $P_{-2419}^{\gamma_2}\in\Fp[X]$, the 
cube $x^3\in\Fp$ is the $j$-invariant of a curve with $N$ points.

The value of the double $\eta$-quotient
$f = {\eta(z/5) \eta (z/13) \over \eta(z) \eta(z/65)}$ at
$z = {-21+\sqrt{-2419}\over2}$ generates the Hilbert class field $H_{-2419}$. 
The minimal polynomial $\Psi$ of $f$ over $\C(j)$ can be computed as in
[\ESB]. It has degree 4 in $j$ and 
degree $84$ in $X$, and we have $r(f)=84/4=21$.
Indeed, the polynomial
$$
P_{-2419}^{f}=
 X^8 + 87 X^7 + 14637 X^6-3810 X^5 + 39662 X^4 + 42026 X^3 + 12593 X^2
-221 X + 1 
$$
has coefficients of no more than $119/21<6$ digits, and its roots
generate $H_{-2419}$ over $\Q(\sqrt{-2419})$. This polynomial splits 
completely modulo $p$. Let $\alpha\in\Fp$ be a root. The polynomial
$\Psi(\alpha,X)\in\Fp[X]$ has degree $4$, and one of its roots in $\Fp$ 
is the $j$-invariant of a curve with $N$ points. 

If we are only interested in an elliptic curve whose group order is a prime
of $60$ decimal digits, we can also use the naive algorithm of trying
random curves over a field $\Fp$ with $p \approx 10^{60}$. For
$p = 10^{60} + 7 = \text{\rm nextprime}(10^{60})$, the smallest positive
integer $j$ such that $j \bmod p$ is the $j$-invariant of a curve of
prime order is $j=180$.

Alternatively, we can use Algorithm 3.1 to construct an elliptic curve 
with CM by $\Z[\zeta_3]$ that has prime order of the desired size. 
In Step 2 we consider consecutive primes $p = 10^{60}+99,\ldots$ congruent
to $1 \bmod 3$.
The fourth prime, $p = 10^{60} + 1059$, yields the prime value
$$
q = 999999999999999999999999999998130705774503095542609960125197,
$$
and the elliptic curve defined by $Y^2 = X^3 + 537824$ has $q$ rational points
over $\Fp$.\par
\ \par\noindent
There is some concern that elliptic curves with `small' endomorphism ring
might be less secure for cryptographic purposes. It is recommended that
the class number of the endomorphism ring should be at least $200$. This is
no problem for our algorithms. Indeed, we can stop in Step 3b in Algorithm 2.2 
only if the class number $h(D)=\deg P_D$ of the discriminant $D$ obtained is
large enough.
In our example, the next higher values after $D=-2419$ for which
we also find solutions to $x^2-Dy^2 = 4N$ with $N+1\pm x$ 
prime are $D =-21003 = -3\cdot 7001, D=-517147=-587\cdot 881$ and $D = 
-590971 = -17\cdot 34763$. The class group of the order of discriminant
$-590971$ is cyclic of order $228>200$.

The primes $N$ needed for cryptography are rather `small' as input for the
two Algorithms in Section 2. 
There is no problem in feeding $N = \text{\rm nextprime}
(10^{2006}) = 10^{2006}+2247$ to Algorithm 2.2.
It yields the prime discriminant $D = -15907$ of class number 15.
The corresponding class polynomial $P_D$ has coefficients of up to 273 digits,
and the desired elliptic curve is readily found.
We do not print it here, as it has coefficients modulo a prime of 2006 digits
that are not particularly pleasing to the human eye.

\enddocument